\documentclass{amsart}
\usepackage{amssymb}

\usepackage{amscd}
\usepackage{thmdefs}

\input{tcilatex}
\theoremstyle{definition}
\theoremstyle{remark}
\numberwithin{equation}{section}

\input tcilatex

\begin{document}
\title[Graph Free Probability Spaces]{Graph Free Product of Noncommutative Probability Spaces}
\author{Ilwoo Cho}
\address{Saint Ambrose Univ., Dep. of Math, Davenport, Iowa, U. S. A}
\email{choillwoo@sau.edu}
\thanks{I appreciate all support from Saint Ambrose University.}
\date{}
\dedicatory{}
\thanks{}
\maketitle

\begin{abstract}
In this paper, we will introduce a new operator-algebraic probability
structure, so-called the Graph Free Product Spaces. Let $G$ be a simplicial
finite graph with its probability-space-vertices $\{(A_{v},$ $\varphi _{v})$ 
$:$ $v$ $\in $ $V(G)\}.$ Define the graph free product $A$ $=$ $\underset{%
v\in V(G)}{*^{G}}$ $A_{v},$ as a vector space $\Bbb{C}$ $\oplus $ $(%
\underset{w\in \Bbb{F}^{+}(G)}{\oplus }$ $A_{w}),$ where $\Bbb{F}^{+}(G)$ is
the free semigroupoid of $G,$ consisting of all vertices as units and all
admissible finite paths, with $A_{w}$ $=$ $A_{w},$ if $w$ is a vertex in $G,$
and with $A_{w}$ $=$ $A_{v_{1}}$ $*$ $A_{v_{2}}$ $*$ ... $*$ $A_{v_{k}},$
whenever $w$ $=$ $[v_{1},$ ..., $v_{k}]$ is a finite path. Also, define the
canonical subalgebra $D^{G}$ of $A$ by $D^{G}$ $=$ $\Bbb{C}$ $\oplus $ $(%
\underset{w\in \Bbb{F}^{+}(G)}{\oplus }$ $\Bbb{C}).$ Then definitely $D^{G}$
is a subalgebra of $A.$ The algebraic pair $(A,$ $E)$ is a noncommutative
probability space with amalgamation over $D^{G}$ in the sense of Voiculescu,
where $E$ $=$ $\underset{w\in \Bbb{F}^{+}(G)}{\oplus }$ $\varphi _{w}$ is
the conditional expectation from $A$ onto $D^{G}.$ In fact, this structure
is a direct producted noncommutative probability space introduced in [9].
This amalgamated noncommutative probability space $(A,$ $E)$ is called the
graph free product space of $\{(A_{v},$ $\varphi _{v})$ $:$ $v$ $\in $ $%
V(G)\}.$ We will consider the noncommutative probability on it. We can
characterize the graph-freeness pictorially for the given graph. i.e., the
subalgebras $A_{w_{1}}$ and $A_{w_{2}}$ are graph-free in $(A,$ $E)$ if and
only if $w_{1}$ and $w_{2}$ are disjoint on the graph $G.$ By using this
graph-freeness, we establish the graph R-transform calculus.
\end{abstract}

\strut \strut

In this paper, we will define and observe Graph Free Probability Spaces
which are the graph free product of probability-space-vertices. Graph free
product spaces are the direct sum of free products of noncommutative
probability spaces, where the direct sum is highly depending on the given
graphs. Roughly speaking, a graph free product space is a certain
graph-depending direct sum of free product spaces. Free Probability has been
developed by Voiculescu, Speicher and various mathematicians from 1980's.
There are two approaches to study it. One of them is the original
Voiculescu's pure analytic approach (See [4]) and the other one is the
Speicher and Nica's combinatorial approach (See [1], [11] and [12]). Let $B$ 
$\subset $ $A$ be unital algebras with $1_{A}=1_{B}$. Suppose that there
exists a conditional expectation $E$ $:$ $A$ $\rightarrow $ $B$ satisfying
the bimodule map property and

\strut

(i) \ \ $E(b)=b,$ for all $b\in B$

(ii) \ $E(bab^{\prime })=bE(a)b^{\prime },$ for all $b,b^{\prime }\in B$ and 
$a\in A.$

\strut

Then the algebraic pair $(A,E)$ is called a noncommutative probability space
with amalgamation over $B$ (or an amalgamated noncommutative probability
space over $B.$ See [11]). All elements in $(A,$ $E)$ are said to be $B$%
-valued random variables. When $B=\Bbb{C}$ and $E$ is a linear functional,
then we call this structure a (scalar-valued) noncommutative probability
space (in short, a probability space) and the elements of $A$ are said to be
(free) random variables. Let $a$ $\in $ $(A,$ $E)$ be a $B$-valued random
variable. Then it contains the following (equivalent) free probabilistic
data,

\strut

\begin{center}
$E(b_{1}a...b_{n}a)$
\end{center}

and

\begin{center}
$k_{n}(b_{1}a,...,b_{n}a)\overset{def}{=}\underset{\pi \in NC(n)}{\sum }%
E_{\pi }(b_{1}a,...,b_{n}a)\mu (\pi ,1_{n}),$
\end{center}

\strut

which are called the $B$-valued $n$-th moment of $a$ and the $B$-valued $n$%
-th cumulant of $a,$ respectively, for all $n$ $\in $ $\Bbb{N}$ and $b_{1},$
..., $b_{n}$ $\in B$ are arbitrary, where $E_{\pi }$ $(...)$ is the
partition-dependent $B$-valued moment of $a$ and $NC(n)$ is the collection
of all noncrossing partitions over $\{1,$ ..., $n\}$ and $\mu $ is the M\"{o}%
bius functional in the incidence algebra $I_{2},$ as the convolution inverse
of the zeta functional $\zeta $,

\strut

\begin{center}
$\zeta (\pi _{1},\pi _{2})\overset{def}{=}\left\{ 
\begin{array}{lll}
1 &  & \text{if }\pi _{1}\leq \pi _{2} \\ 
0 &  & \text{otherwise,}
\end{array}
\right. $
\end{center}

\strut

for all $\pi _{1},$ $\pi _{2}\in NC(n),$ and $n\in \Bbb{N}$ (See [11]). When 
$b_{1}$ $=$ ... $=$ $b_{n}$ $=1_{B},$ for $n\in \Bbb{N},$ we say that $%
E(a^{n})$ and $k_{n}(a,$ $...,$ $a)$ are trivial $n$-th $B$-valued moment
and cumulant of $a,$ respectively. Recall that the collection $NC(n)$ is a
lattice with the total ordering $\leq $ (See [1], [2], [10] and [11]).
Again, if $B$ $=$ $\Bbb{C},$ then we have the same (scalar-valued) moments
and cumulants of each random variable $a$ $\in $ $(A,$ $E).$ But, in this
case, since $B$ $=$ $\Bbb{C}$ commutes with $A,$ we have that

\strut

\begin{center}
$E\left( b_{1}a...b_{n}a\right) =\left( b_{1}...b_{n}\right) E(a^{n})$
\end{center}

and

\begin{center}
$k_{n}\left( b_{1}a,...,b_{n}a\right) =\left( b_{1}...b_{n}\right)
k_{n}(a,...,a),$
\end{center}

\strut

for all $n\in \Bbb{N}$ and for all $b_{1},...,b_{n}\in B=\Bbb{C}.$ These
mean that we only need to consider the trivial moments and cumulants of $a,$
with respect to $E$ $:$ $A$ $\rightarrow $ $B$ $=$ $\Bbb{C}.$ So, the
(scalar-valued) moments and cumulants of $a$ are defined by the trivial $B$ $%
=$ $\Bbb{C}$-valued moments and cumulants of it. More generally, if an
arbitrary unital algebra $B$ commutes with $A,$ where $A$ is over $B,$ then
we can verify that the trivial $B$-valued moments and cumulants of a $B$%
-valued random variable $a$ contains the full free probabilistic data of $a.$
So, in this case, without loss of generality, we define the $B$-valued
moments and cumulants of $a$ by the trivial $B$-valued moments and cumulants
of $a.$

\strut \strut

Let $A_{1}$ and $A_{2}$ be subalgebras of $A.$ We say that they are free
over $B$ if all mixed cumulants of $A_{1}$ and $A_{2}$ vanish. Let $S_{1}$
and $S_{2}$ be subsets of $A.$ We say that subsets $S_{1}$ and $S_{2}$ are
free over $B$ if the subalgebras $A_{1}$ $=$ $A\lg $ $\{S_{1},$ $B\}$ and $%
A_{2}$ $=$ $A\lg $ $\{S_{2},$ $B\}$ are free over $B.$ In particular, the $B$%
-valued random variables $a_{1}$ and $a_{2}$ are free over $B$ if the
subsets $\{a_{1}\}$ and $\{a_{2}\}$ are free over $B$ in $(A,$ $E).$
Equivalently, given two $B$-valued random variables $a_{1}$ and $a_{2}$ are
free over $B$ in $(A,$ $E)$ if all mixed $B$-valued cumulants of $a_{1}$ and 
$a_{2}$ vanish. (Recall that, when $A$ is a $*$-algebra, they are free over $%
B$ if all mixed $B$-valued cumulants of $P(a_{1},$ $a_{1}^{*})$ and $%
Q(a_{2}, $ $a_{2}^{*})$ vanish, for all $P,$ $Q$ in $\Bbb{C}$ $[z_{1},$ $%
z_{2}].$)

\strut

Let $G$ be a finite simplicial (undirected) graph with its vertex set $V(G)$
and its edge set $E(G).$ A graph is simplicial if it has neither multiple
edges between two vertices nor loop-edges. For example, the circulant graph $%
C_{N}$ with $N$-vertices is a simplicial graph. Let $e$ be an edge
connecting vertices $v_{1}$ and $v_{2}.$ Then we denote $e$ by $[v_{1},$ $%
v_{2}]$ or $[v_{2},$ $v_{1}],$ and we assume that $[v_{1},$ $v_{2}]$ $=$ $%
[v_{2},$ $v_{1}].$ Suppose $w$ is a finite path in the graph connecting
edges $[v_{1},$ $v_{2}],$ $[v_{2},$ $v_{3}],$ ..., $[v_{k-1},$ $v_{k}].$
Then, for convenience, we denote $w$ by $[v_{1},$ $v_{2},$ ..., $v_{k}].$ Of
course, we assume that

\strut

\begin{center}
$[v_{1},v_{2},...,v_{k}]=[v_{k},...,v_{2},v_{1}],$
\end{center}

\strut

as a same finite path. In this case, the length $\left| w\right| $ of the
finite path $w$ is defined to be $k$ $-$ $1,$ the cardinality of the
admissible edges. Define the set $FP(G)$ by the set of all such finite paths
and define the free semigroupoid $\Bbb{F}^{+}(G)$ of the given graph $G$ by

\strut

\begin{center}
$\Bbb{F}^{+}(G)=V(G)\cup FP(G),$
\end{center}

\strut

as a set.

\strut

Let $G$ be a finite simplicial graph with its probability-space-vertices

\strut

\begin{center}
$\{(A_{v},\varphi _{v}):v\in V(G)\},$
\end{center}

\strut

where $(A_{v},\varphi _{v})$ are noncommutative probability space with their
linear functionals $\varphi _{v}$ $:$ $A_{v}$ $\rightarrow $ $\Bbb{C},$
indexed by the vertices $v$ in $V(G).$ Define the graph free (or $G$-free)
product $A^{G}$ of $A_{v}$'s by

\strut

\begin{center}
$A^{G}\overset{denote}{=}\,\underset{v\in V(G)}{*^{G}}A_{v}\overset{def}{=}%
\Bbb{C}\oplus \left( \underset{w\in \Bbb{F}^{+}(G)}{\oplus }A_{w}\right) ,$
\end{center}

\strut

as an algebra, where $A_{w}=A_{w},$ if $w$ $\in $ $V(G),$ and

\strut

\begin{center}
$A_{w}=A_{[v_{1},...,v_{k}]}=A_{v_{1}}*...*A_{v_{k}},$
\end{center}

\strut

for all $w=[v_{1},...,v_{k}]\in FP(G),$ where $``*"$ is the usual
(scalar-valued) free product. Similarly, define the free product of linear
functional $\varphi _{w}$ on $A_{w}$ by

\strut

\begin{center}
$\varphi _{w}=\varphi _{w}$ if $w\in V(G)$
\end{center}

and

\begin{center}
$\varphi _{w}=\varphi _{[v_{1},...,v_{k}]}=\varphi _{v_{1}}*...*\varphi
_{v_{k}},$
\end{center}

\strut

if $w=[v_{1},...,v_{k}]\in FP(G).$ Then, for each $w$ $\in $ $\Bbb{F}%
^{+}(G), $ we have the free product space $(A_{w},$ $\varphi _{w}),$ as a
noncommutative probability space with its linear functional $\varphi _{w}.$
Notice that if $w_{1}$ and $w_{2}$ are finite paths in $FP(G)$ and assume
that $w_{3}$ $=$ $w_{1}w_{2}$ is also an admissible finite path (i.e., $%
w_{3} $ $\in $ $FP(G)$). Then we have that

\strut

(0.1) $\ \ \ \ \ \ \ \ \ \ \ \left( A_{w_{3}},\varphi _{w_{3}}\right)
=\left( A_{w_{1}},\varphi _{w_{1}}\right) *\left( A_{w_{2}},\varphi
_{w_{2}}\right) .$

\strut

On the other hands, if $w_{1}$ and $w_{2}$ are finite paths in $FP(G)$ and
if they are not admissible (i.e., $w_{1}w_{2}$ $\notin $ $FP(G)$ or $%
w_{2}w_{1}$ $\notin $ $FP(G)$). Then the free product $A_{w_{1}}$ $*$ $%
A_{w_{2}}$ of $A_{w_{1}}$ and $A_{w_{2}}$ does not exist as a direct summand
of the graph free product $A^{G}$ $=$ $\underset{v\in V(G)}{*^{G}}$ $A_{v},$
however, there exists a direct summand $A_{w_{1}}$ $\oplus $ $A_{w_{2}}$ of $%
A^{G}.$

\strut

Define the following canonical subalgebra $D^{G}$ of $A^{G}$ by

\strut

\begin{center}
$D^{G}\overset{def}{=}\Bbb{C}\oplus \left( \underset{w\in \Bbb{F}^{+}(G)}{%
\oplus }\Bbb{C}_{w}\right) ,$ \ with \ $\Bbb{C}_{w}=\Bbb{C},$ $\forall w\in 
\Bbb{F}^{+}(G).$
\end{center}

\strut

Then we can define the conditional expectation $E$ from $A^{G}$ onto $D^{G}$
by $E$ $=$ $\underset{w\in \Bbb{F}^{+}(G)}{\oplus }$ $\varphi _{w}.$ The
algebraic pair $(A^{G},$ $E)$ is called the graph free product of
probability-space vertices $\{(A_{v},$ $\varphi _{v})$ $:$ $v$ $\in $ $%
V(G)\} $ or the $G$-free probability space. And all elements of this $G$%
-free product space $(A^{G},$ $E)$ are called $G$-free random variables.

\strut

The main purpose of studying the $G$-free probability spaces is to extend
the classical free probabilistic data to a certain free-product-like
structure depending on a pure combinatorial object (a graph). Also, the $G$%
-free product, itself, is the extension of the usual free product. For
example, if the given graph $H$ is complete, in the sense that there always
exists an edge $[v,$ $v^{\prime }],$ for all pair $(v,$ $v^{\prime })$ of
vertices, then the $H$-free product is nothing but the classical free
product.

\strut

In Chapter 1, we will review the direct producted noncommutative probability
spaces introduced in [9]. In [9], we only observed the finite direct product
of noncommutative probability spaces, but in Chapter 1, we will consider the
infinite direct product of noncommutative probability spaces. After
considering the direct producted noncommutative probability spaces, we can
realize that our graph free product $(A^{G},$ $E)$ of
probability-space-vertices $\{(A_{v},$ $\varphi _{v})\}_{v\in V(G)}$ is a
certain countable direct producted noncommutative probability space, having
its direct summands as free products of probability spaces, depending on an
admissible finite path.

\strut

In Chapter 2, we consider the graph-freeness characterization. And we can
characterize the graph-freeness of two subalgebras of $A^{G},$ in terms of
the disjointness on the graph $G.$ i.e., the subalgebras $A_{w_{1}}$ and $%
A_{w_{2}}$ of $A^{G}$ are graph-free if and only if $w_{1}$ and $w_{2}$ are
disjoint on the graph $G.$ Also, we will define the graph moments and graph
cumulants of a $G$-free random variable in $(A^{G},$ $E).$ And then, we will
consider the $G$-freeness on $A^{G}.$ Under this $G$-freeness, we will do
the graph R-transform calculus over $D^{G}.$ In Chapter 3, we will study
certain $G$-free random variables under the $W^{*}$-setting. Finally, in
Chapter 4, we will consider certain subalgebras of the graph free
probability space.

\strut \strut

\strut \strut

\strut

\section{Direct Producted Noncommutative Probability Spaces}

\strut

\strut

Throughout this chapter, let's fix a sufficiently big number $N\in \Bbb{N}$
(possibly, $N$ $\rightarrow $ $\infty $) and the collection $\mathcal{F}$ of
(scalar-valued) noncommutative probability spaces,

\strut

\begin{center}
$\mathcal{F}=\{(A_{i},\varphi _{i}):i=1,...,N\}.$
\end{center}

\strut

Then, for the given unital algebras $A_{1},$ ..., $A_{N}$ in $\mathcal{F},$
we can define the direct producted unital algebra

$\strut $

\begin{center}
$A$ $=$ $\oplus _{j=1}^{N}$ $A_{j}=\{\oplus _{j=1}^{N}a_{j}:a_{j}\in A_{j},$ 
$j=1,...,N\},$
\end{center}

\strut

where $\oplus $ means the algebra direct sum with its componentwise vector
addition and the vector multiplication defined again componentwisely by

\strut

\begin{center}
$\left( \oplus _{j=1}^{N}a_{j}\right) \cdot \left( \oplus
_{j=1}^{N}a_{j}^{\prime }\right) =\oplus _{j=1}^{N}a_{j}a_{j}^{\prime },$
\end{center}

\strut

for all $\oplus _{j=1}^{N}a_{j},$ $\oplus _{j=1}^{N}a_{j}^{\prime }$ $\in $ $%
A.$ This direct product $A$ of $A_{1},$ ..., $A_{N}$ is again a unital
algebra with its unity $1_{A}=\oplus _{j=1}^{N}1.$

\strut \strut

Define the subalgebra $D_{N}$ of the direct producted algebra $A$ $=$ $%
\oplus _{j=1}^{N}$ $A_{j}$ by

\strut

\begin{center}
$D_{N}\overset{def}{=}\oplus _{j=1}^{N}\Bbb{C}_{j},$ with $\Bbb{C}_{j}=\Bbb{C%
},$ $\forall j=1,...,N.$
\end{center}

\strut

Then $D_{N}$ is a commutative subalgebra of $M_{N}(\Bbb{C})$ and it is
isomorphic to the matricial algebra $\Delta _{N}$ generated by all $N$ $%
\times $ $N$ diagonal matrices in the matricial algebra $M_{N}(\Bbb{C}).$ We
will call this subalgebra $D_{N}$ of $A,$ the $N$-th diagonal algebra of $A.$

\strut

\begin{definition}
Define the direct producted noncommutative probability space $A$ of unital
algebras $A_{1},$ ..., $A_{N},$ by the noncommutative probability space $(A,$
$E)$ with amalgamation over the $N$-th diagonal algebra $D_{N},$ where $A$ $=
$ $\oplus _{j=1}^{N}$ $A_{j}$ is the direct product of $A_{1},$ ..., $A_{N}$
and $E$ $:$ $A$ $\rightarrow $ $D_{N}$ is the conditional expectation from $A
$ onto $D_{N}$ defined by

\strut 

$\ \ \ \ \ \ \ \ \ \ \ \ \ \ \ \ \ \ E\left( \oplus _{j=1}^{N}a_{j}\right)
=\oplus _{j=1}^{N}\varphi _{j}(a_{j}),$

\strut 

for all $\oplus _{j=1}^{N}a_{j}\in A.$ Some times, we will denote $E$ by $%
\oplus _{j=1}^{N}$ $\varphi _{j}.$
\end{definition}

\strut

It is easy to see that the $\Bbb{C}$-linear map $E$ is indeed a conditional
expectation;

\strut

(i) \ \ $E\left( \oplus _{j=1}^{N}\alpha _{j}\right) =\oplus
_{j=1}^{N}\varphi _{j}(\alpha _{j})=\oplus _{j=1}^{N}\alpha _{j},$

\strut

\ \ \ \ \ \ for all $\oplus _{j=1}^{N}\alpha _{j}\in D_{N}.$

\strut

(ii) \ $E\left( (\oplus _{j=1}^{N}\alpha _{j})\left( \oplus
_{j=1}^{N}a_{j}\right) (\oplus _{j=1}^{N}\alpha _{j}^{\prime })\right) $

\strut

$\ \ \ \ \ \ \ =E\left( \oplus _{j=1}^{N}\alpha _{j}a_{j}\alpha _{j}^{\prime
}\right) =\oplus _{j=1}^{N}\varphi _{j}(\alpha _{j}a_{j}\alpha _{j}^{\prime
})$

\strut

$\ \ \ \ \ \ =\oplus _{j=1}^{N}\left( \alpha _{j}\cdot \varphi
_{j}(a_{j})\cdot \alpha _{j}^{\prime }\right) $

\strut

$\ \ \ \ \ \ \ =\left( \oplus _{j=1}^{N}\alpha _{j}\right) \cdot \left(
\oplus _{j=1}^{N}\varphi _{j}(a_{j})\right) \cdot \left( \oplus
_{j=1}^{N}\alpha _{j}^{\prime }\right) $

\strut

$\ \ \ \ \ \ \ =(\oplus _{j=1}^{N}\alpha _{j})\left( E(\oplus
_{j=1}^{N}a_{j}))\right) (\oplus _{j=1}^{N}\alpha _{j}^{\prime }),$

\strut

\ \ \ \ \ \ for all $\oplus _{j=1}^{N}\alpha _{j},$ $\oplus _{j=1}^{N}\alpha
_{j}^{\prime }$ $\in $ $D_{N}$ and $\oplus _{j=1}^{N}a_{j}\in A.$

\strut

By (i) and (ii), the map $E$ is a conditional expectation from $A$ $=$ $%
\oplus _{j=1}^{N}$ $A_{j}$ onto $D_{N}.$ Thus the algebraic pair $(A,$ $E)$
is a noncommutative probability space with amalgamation over the $N$-th
diagonal algebra $D_{N}.$

\strut

\begin{definition}
Let $(A_{j},\varphi _{j})$ be $*$-probability spaces, for $j$ $=$ $1,$ ..., $%
N,$ where $A_{j}$'s are unital $*$-algebras and $\varphi _{j}$'s are linear
functional satisfying that $\varphi _{j}(a_{j}^{*})$ $=$ $\overline{\varphi
_{j}(a_{j})},$ for all $a_{j}$ $\in $ $A_{j},$ for $j$ $=$ $1,$ ..., $N.$
Then the direct producted $*$-probability space $(A,$ $E)$ is defined by the
algebraic pair of direct product $A$ $=$ $\oplus _{j=1}^{N}$ $A_{j}$ of $%
A_{1},$ ..., $A_{N},$ as a $*$-algebra, with $(\oplus _{j=1}^{N}a_{j})^{*}$ $%
=$ $\oplus _{j=1}^{N}a_{j}^{*}$, for $\oplus _{j=1}^{N}a_{j}$ $\in $ $A,$
and the conditional expectation satisfying the above condition (i), (ii) and
the following condition (iii);

\strut 

$\ \ \ \ \ \ \ \ \ \ E\left( (\oplus _{j=1}^{N}a_{j})^{*}\right) =E\left(
\oplus _{j=1}^{N}a_{j}\right) ^{*}$ \ in $D_{N},$

\strut 

for all $\oplus _{j=1}^{N}a_{j}\in A.$ Similarly, we can define the direct
producted $C^{*}$-probability (or $W^{*}$-probability) spaces, if $A$ is a $%
C^{*}$-direct sum (resp. $W^{*}$-direct sum) of $C^{*}$-algebras (resp. von
Neumann algebras), also denoted by $\oplus _{j=1}^{N}$ $A_{j},$ and the
conditional expectation $E$ satisfy (i), (ii) and (iii), with the continuity
under the given $C^{*}$-topology (resp. $W^{*}$-topology) on $A$. (Notice
that the $C^{*}$ or $W^{*}$-topology of $A$ is the product topology of those
of $A_{1},$ ..., $A_{N}$ ($N$ $\rightarrow $ $\infty $). And the continuity
comes from that of $\varphi _{1},$ ..., $\varphi _{N},$ under the product
topology ($N$ $\rightarrow $ $\infty $).)
\end{definition}

\strut \strut

In the rest of this section, we will let $A_{1},$ ..., $A_{N}$ be just
unital algebras without the involution and topology. However, if we put
involution or topology on $A_{1},$ ..., $A_{N},$ we would have the same or
similar results.

\strut

Now, we will consider the $D_{N}$-freeness on the direct producted
noncommutative probability space $(A,$ $E)$. Notice that the $N$-th diagonal
algebra $D_{N}$ satisfies that

\strut

(1.1) $\ \ \ \ \ \ \ \ \ dx=xd,$ for all $d\in D_{N}$ and $x\in A$,

\strut

as a subalgebra of our direct product $A$ $=$ $\oplus _{j=1}^{N}$ $A_{j}.$
By (1.1) and the commutativity of $D_{N},$ we only need to consider the
trivial $D_{N}$-valued moments and cumulants of $D_{N}$-valued random
variables, for studying the free probabilistic data of them. i.e., we have
that, for any $a\in A,$

\strut

(1.2) $\ \ \ \ \ \ \ \ \ \ \ \ \ E(d_{1}a...d_{n}a)=\left(
d_{1}...d_{n}\right) E(a^{n})$

\strut

and

\strut

(1.3) \ \ \ \ \ \ \ \ $k_{n}\left( d_{1}a,...,d_{n}a\right) =\left(
d_{1}...d_{n}\right) k_{n}\left( a,...,a\right) ,$

\strut

for all $n\in \Bbb{N}$ and for any arbitrary $d_{1},...,d_{n}\in D_{N}.$ So,
the relations (1.2) and (1.3) shows that it is enough to consider the
trivial $D_{N}$-valued moments $E(a^{n})$ and cumulants $k_{n}(a,$ ..., $a)$
of $D_{N}$-valued random variables ($n$ $\in $ $\Bbb{N}$), whenever we want
to know about the free probabilistic information of those $D_{N}$-valued
random variables.

\strut

\begin{proposition}
Let $(A,E)$ be the direct producted noncommutative probability space with
amalgamation over the $N$-th diagonal algebra $D_{N},$ where $A$ $=$ $\oplus
_{j=1}^{N}$ $A_{j}$ and let $a_{1}$ and $a_{2}$ be $D_{N}$-valued random
variables in $(A,$ $E).$ Then they are free over $D_{N}$ in $(A,$ $E)$ if
and only if all mixed \textbf{trivial }$D_{N}$-valued cumulants of them
vanish. $\square $
\end{proposition}

\strut

\begin{remark}
If we have a direct producted $*$-probability (or $C^{*}$-probability or $%
W^{*}$-probability) space $(A,$ $E),$ then the above proposition does not
hold, in general. But, we can get that the $D_{N}$-valued random variables $%
a_{1}$ and $a_{2}$ are free over $D_{N}$ in $(A,$ $E)$ if and only if all
mixed trivial $D_{N}$-valued cumulants of $a_{1},$ $a_{1}^{*},$ $a_{2}$ and $%
a_{2}^{*}$ vanish. $\square $
\end{remark}

\strut \strut

We will consider the $D_{N}$-valued moments of an arbitrary random variable
in the direct producted noncommutative probability space.

\strut

\begin{proposition}
Let $\left( A=\oplus _{j=1}^{N}A_{j},\text{ }E=\oplus _{j=1}^{N}\varphi
_{j}\right) $ be the direct producted noncommutative probability space over
the $N$-th diagonal algebra $D_{N}$ and let $x$ $=$ $\oplus _{j=1}^{N}$ $%
a_{j}$ be the $D_{N}$-valued random variable in $(A,$ $E).$ Then the trivial 
$n$-th moment of $x$ is

\strut 

(1.4)$\ \ \ \ \ \ \ \ \ \ \ \ \ \ \ \ \ \ E(x^{n})=\oplus _{j=1}^{N}\left(
\varphi _{j}(a_{j}^{n})\right) ,$

\strut 

for all $n\in \Bbb{N}.$ \ $\square $
\end{proposition}

\strut

The above proposition is easily proved by the fact that

\strut

\begin{center}
$x^{n}=\left( \oplus _{j=1}^{N}a_{j}\right) ^{n}=\oplus _{j=1}^{N}a_{j}^{n},$
for all $n\in \Bbb{N}.$
\end{center}

\strut

It shows that if we know the $n$-th moments of $a_{j}$ $\in $ $(A_{j},$ $%
\varphi _{j}),$ for $j$ $=$ $1,$ ..., $N,$ then we can compute the $D_{N}$%
-valued moment of the $D_{N}$-valued random variable $\oplus _{j=1}^{N}$ $%
a_{j}$ in $(A,$ $E).$ In \ other words, the free probabilistic information
of $A_{1}$, ..., $A_{N}$ affects the amalgamated free probabilistic
information of $\oplus _{j=1}^{N}$ $A_{j}.$ Now, let's compute the trivial $%
n $-th cumulant of an arbitrary $D_{N}$-valued random variable;

\strut

\begin{proposition}
Let $\left( A=\oplus _{j=1}^{N}A_{j},E=\oplus _{j=1}^{N}\varphi _{j}\right) $
be the direct producted noncommutative probability space over the $N$-th
diagonal algebra $D_{N}$ and let $x$ $=$ $\oplus _{j=1}^{N}a_{j}$ be the $%
D_{N}$-valued random variable in $(A,$ $E).$ Then the trivial $n$-th
cumulant of $x$ is

\strut \strut 

\ \ \ \ \ $\ \ \ k_{n}\left( \underset{n\text{-times}}{\underbrace{%
x,........,x}}\right) =\oplus _{j=1}^{N}\left( k_{n}^{(j)}\left(
a_{j},...,a_{j}\right) \right) ,$

\strut 

for all $n\in \Bbb{N},$ where $k_{n}^{(i)}(...)$ is the $n$-th cumulant
functional with respect to the noncommutative probability space $(A_{i},$ $%
\varphi _{i}),$ for all $i$ $=$ $1,$ ..., $N.$
\end{proposition}

\strut

\begin{proof}
Fix $n\in \Bbb{N}.$ Then

\strut

$\ \ k_{n}\left( x,...,x\right) =\underset{\pi \in NC(n)}{\sum }E_{\pi
}(x,...,x)\mu (\pi ,1_{n})$

\strut

(1.5)\ \ \ \ 

\strut

$\ \ \ \ \ \ \ =\underset{\pi \in NC(n)}{\sum }\left( \underset{V\in \pi }{%
\Pi }E_{V}(x,...,x)\right) \mu (\pi ,1_{n})$

\strut

by (1,1), where $E_{V}(x,...,x)=E(x^{\left| V\right| }),$ where $\left|
V\right| $ is the length of the block (See [1] and [11])

\strut

$\ \ \ \ \ \ \ =\underset{\pi \in NC(n)}{\sum }\left( \underset{V\in \pi }{%
\Pi }\left( \oplus _{j=1}^{N}(\varphi _{j}(a_{j}^{\left| V\right| }))\right)
\right) \mu (\pi ,1_{n})$

\strut

by (1.4)

\strut

$\ \ \ \ \ \ \ \ =\underset{\pi \in NC(n)}{\sum }\left( \oplus
_{j=1}^{N}\left( \underset{V\in \pi }{\Pi }\varphi _{j}(a_{j}^{\left|
V\right| })\right) \right) \mu (\pi ,1_{n})$

\strut

since $\left( \oplus _{j=1}^{N}\alpha _{j}\right) \cdot \left( \oplus
_{j=1}^{N}\alpha _{j}^{\prime }\right) =\oplus _{j=1}^{N}\alpha _{j}\alpha
_{j}^{\prime }$ in $D_{N}$

\strut

$\ \ \ \ \ \ \ \ =\underset{\pi \in NC(n)}{\sum }\left( \oplus
_{j=1}^{N}\left( (\underset{V\in \pi }{\Pi }\varphi _{j}(a_{j}^{\left|
V\right| }))\mu (\pi ,1_{n})\right) \right) ,$

\strut

since the $N$-th diagonal algebra $D_{N}$ is a vector space (i.e., if we let 
$\mu _{\pi }$ $=$ $\mu (\pi ,$ $1_{n})$ in $\Bbb{C},$ for the fixed $\pi $ $%
\in $ $NC(n),$ then $(\oplus _{j=1}^{N}$ $\alpha _{j})$ $\cdot $ $\mu _{\pi
} $ $=$ $\oplus _{j=1}^{N}$ $(\alpha _{j}\mu _{\pi }),$ for all $\oplus
_{j=1}^{N}$ $\alpha _{j}$ $\in $ $D_{N}.$)

\strut

$\ \ \ \ \ \ \ \ =\oplus _{j=1}^{N}\left( \underset{\pi \in NC(n)}{\sum }%
\left( \underset{V\in \pi }{\Pi }\varphi _{j}(a_{j}^{\left| V\right|
})\right) \mu (\pi ,1_{n})\right) $

\strut

since $(\alpha _{1},$ ..., $\alpha _{N})$ $+$ $(\alpha _{1}^{\prime },$ ..., 
$\alpha _{N}^{\prime })$ $=$ $(\alpha _{1}$ $+$ $\alpha _{1}^{\prime },$
..., $\alpha _{N}$ $+$ $\alpha _{N}^{\prime })$ in $D_{N}.$

\strut

\strut (1.6)

\strut

$\ \ \ \ \ \ \ \ =\oplus _{j=1}^{N}\left(
k_{n}^{(j)}(a_{j},...,a_{j})\right) ,$

\strut

where $k_{n}^{(i)}(...)$ is the (scalar-valued) $n$-th cumulant functional
with respect to the (scalar-valued) noncommutative probability space $%
(A_{i}, $ $\varphi _{i}),$ for all $i$ $=$ $1,$ ..., $N.$
\end{proof}

\strut

Therefore, the equality (1.6) shows that the $n$-th $D_{N}$-valued cumulant

$\strut $

\begin{center}
$k_{n}\left( \oplus _{j=1}^{N}a_{j},...,\oplus _{j=1}^{N}a_{j}\right) $
\end{center}

\strut

of the $D_{N}$-valued random variable $\oplus _{j=1}^{N}a_{j}$ in the direct
producted noncommutative probability space $(A,$ $E),$ is nothing but the $N$%
-tuple of $n$-th (scalar-valued) cumulants of $a_{1},$ ..., $a_{N},$

\strut

\begin{center}
$\oplus _{j=1}^{N}\left( k_{n}^{(j)}(a_{j},...,a_{j})\right) $
\end{center}

\strut

in $D_{N}.$ By (1.6), we can get the following $D_{N}$-freeness
characterization on the direct producted noncommutative probability space

$\strut $

\begin{center}
$\left( A=\oplus _{j=1}^{N}A_{j},\text{ }E=\oplus _{j=1}^{N}\varphi
_{j}\right) .$
\end{center}

\strut

\begin{theorem}
Let $\left( A=\oplus _{j=1}^{N}A_{j},E=\oplus _{j=1}^{N}\varphi _{j}\right) $
be the given direct producted noncommutative probability space over the $N$%
-th diagonal algebra $D_{N}$ and let $x_{1}$ $=$ $\oplus _{j=1}^{N}a_{j}$
and $x_{2}$ $=$ $\oplus _{j=1}^{N}b_{j}$ be the $D_{N}$-valued random
variables in $(A,E).$ Then $x_{1}$ and $x_{2}$ are free over $D_{N}$ in $(A,$
$E)$ if and only if $a_{j}$ and $b_{j}$ are free in $(A_{j}$, $\varphi _{j}),
$ for all $j$ $=$ $1,$ ..., $N.$
\end{theorem}

\strut

\begin{proof}
($\Leftarrow $) Assume that random variables $a_{j}$ and $b_{j}$ are free in 
$(A_{j},\varphi _{j}),$ for all $j$ $=$ $1,$ ..., $N.$ Then, by the
freeness, all mixed $n$-th cumulants of $a_{j}$ and $b_{j}$ vanish, for all $%
j$ $=$ $1,$ ..., $N$ and for all $n$ $\in $ $\Bbb{N}$ $\setminus $ $\{1\}.$
By the previous theorem, it is sufficient to show that the $N$-tuples $x_{1}$
$=$ $\oplus _{j=1}^{N}a_{j}$ and $x_{2}$ $=$ $\oplus _{j=1}^{N}$ $b_{j}$ in $%
(A,$ $E)$ have vanishing mixed \textbf{trivial} $D_{N}$-valued cumulants.
Fix $n$ $\in $ $\Bbb{N}$ $\setminus $ $\{1\}$ and let $(x_{i_{1}},$ ..., $%
x_{i_{n}})$ are mixed $n$-tuple of $x_{1}$ and $x_{2},$ where $(i_{1},$ ..., 
$i_{n})$ $\in $ $\{1,$ $2\}^{n}.$ Then

\strut

$\ k_{n}$\strut $\left( x_{i_{1}},...,x_{i_{n}}\right) =k_{n}\left( \oplus
_{j=1}^{N}a_{j}^{i_{1}},...,\oplus _{j=1}^{N}a_{j}^{i_{n}}\right) $

\strut

$\ \ \ \ \ \ \ \ \ \ \ \ \ \ \ =\oplus _{j=1}^{N}\left(
k_{n}^{(j)}(a_{j}^{i_{1}},...,a_{j}^{i_{n}})\right) $

\strut

by the previous proposition

\strut

$\ \ \ \ \ \ \ \ \ \ \ \ \ \ \ =0_{D_{N}},$

\strut

by the hypothesis.

\strut

($\Rightarrow $) Let's assume that the $D_{N}$-valued random variables $%
x_{1} $ $=$ $\oplus _{j=1}^{N}$ $a_{j}$ and $x_{2}$ $=$ $\oplus _{j=1}^{N}$ $%
b_{j}$ are free over $D_{N}$ in $(A,$ $E)$ and assume also that there exists 
$j$ in $\{1,$ ..., $N\}$ such that $a_{j}$ and $b_{j}$ are not free in $%
(A_{j},$ $\varphi _{j}).$ Since $a_{j}$ and $b_{j}$ are not free in $(A_{j},$
$\varphi _{j}),$ there exists $n$ $\in $ $\Bbb{N}$ and the mixed $n$-tuple $%
(i_{1},$ ..., $i_{n})$ $\in $ $\{1,$ $2\}^{n}$ such that $%
k_{n}^{(j)}(a_{j}^{i_{1}},$ ..., $a_{j}^{i_{n}})$ $=$ $\beta $ $\neq $ $0.$
Now, fix the number $n$ and the $n$-tuple $(i_{1},$ ..., $i_{n}).$ Consider
the following mixed trivial $D_{N}$-valued cumulants of $x_{1}$ and $x_{2}$;

\strut

$\ k_{n}\left( x_{i_{1}},...,x_{i_{n}}\right) =k_{n}\left( \oplus
_{j=1}^{N}a_{j}^{i_{1}},...,\oplus _{j=1}^{N}a_{j}^{i_{n}}\right) $

\strut

$\ \ \ \ \ \ \ \ \ \ \ =\oplus _{j=1}^{N}\left(
k_{n}^{(j)}(a_{j}^{i_{1}},...,a_{j}^{i_{n}})\right) $

\strut

$\ \ \ \ \ \ \ \ \ \ \ =0\oplus ...\oplus 0$ $\oplus \ \underset{j\text{-th}%
}{\frame{$\beta $}}\ \oplus 0\oplus ...\oplus 0\neq 0_{D_{N}}.$

\strut

Therefore, there exists the nonvanishing $D_{N}$-valued cumulant of $x_{1}$
and $x_{2}.$ This contradict our assumption that $D_{N}$-valued random
variables $x_{1}$ and $x_{2}$ are free over $D_{N}$ in $(A,$ $E).$
\end{proof}

\strut

\begin{remark}
\strut Now, assume that above $(A,$ $E)$ is direct producted $*$-probability
(or $C^{*}$-probability or $W^{*}$-probability) space of $(A_{1},$ $\varphi
_{1}),$ ..., $(A_{N},$ $\varphi _{N}),$ where $(A_{j},$ $\varphi _{j})$'s
are $*$-probability spaces, for all $j$ $=$ $1,$ ..., $N.$ Then the same
results holds true. But the proof should be different, by the previous
remark. Assume that $x_{1}$ $=$ $\oplus _{j=1}^{N}$ $a_{j}$ and $x_{2}$ $=$ $%
\oplus _{j=1}^{N}$ $b_{j}$ in $(A,$ $E)$ are free over $D_{N}$ in $(A,$ $E).$
Then $*$-Alg$(\{x_{1}\},$ $D_{N})$ and $*$-Alg$(\{x_{2}\},$ $D_{N})$ are
free over $D_{N},$ where $*$-Alg $(S_{1},$ $S_{2})$ means the $*$-algebra
generated by sets $S_{1}$ and $S_{2}.$ Notice that

\strut 

$\ \ \ \ \ \ \ \ \ *$-Alg $(\{x_{1}\},D_{N})=D_{N}\oplus \left( \oplus
_{j=1}^{N}\left( *\text{-Alg}(\{a_{j}\})\right) \right) $

and

$\ \ \ \ \ \ \ \ \ *$-Alg$(\{x_{2}\},D_{N})=D_{N}\oplus \left( \oplus
_{j=1}^{N}\left( *\text{-Alg}(\{b_{j}\})\right) \right) ,$

\strut 

as $*$-subalgebras in $A.$ Assume now that $a_{j_{0}}$ and $b_{j_{0}}$ are
not free in $(A_{j_{0}},$ $\varphi _{j_{0}}),$ for some $j_{0}$ $\in $ $\{1,$
..., $N\}.$ Then the direct summand

\strut 

\ $\ \ \ D_{1}=D_{N}\oplus \left( 0\oplus ...\oplus 0\oplus \,\underset{j_{0}%
\text{-th}}{*\text{-Alg}(\{a_{j_{0}}\})}\,\oplus 0\oplus ...\oplus 0\right) $

and

\ $\ \ \ D_{2}=D_{N}\oplus \left( 0\oplus ...\oplus 0\oplus \,\underset{j_{0}%
\text{-th}}{*\text{-Alg}(\{b_{j_{0}}\})}\,\oplus 0\oplus ...\oplus 0\right) $

\strut 

are not free over $D_{N}$ in $(A,$ $E).$ Since $*$-Alg $(\{x_{i}\},$ $D_{N})$
contain $D_{i},$ for $i$ $=$ $1,$ $2,$ and since $D_{1}$ and $D_{2}$ are not
free over $D_{N},$ $*$-Alg $(\{x_{1}\},$ $D_{N})$ and $*$-Alg $(\{x_{2}\},$ $%
D_{N})$ are not free over $D_{N}.$ This contradict our assumption. The
converse is also similarly proved. Also, we provide the following proof of
the converse;

\strut 

Since $a_{j}$ and $b_{j}$ are free in $(A_{j},\varphi _{j}),$ for all $j$ $=$
$1,$ ..., $N,$ all mixed cumulants of $a_{j}$, $a_{j}^{*},$ $b_{j}$ and $%
b_{j}^{*}$ vanish, for all $j$ $=$ $1,$ ..., $N.$ So, it suffices to show
that all mixed trivial $D_{N}$-valued cumulants of $x_{1},$ $x_{1}^{*}$, $%
x_{2}$ and $x_{2}^{*}$ vanish. Notice that

\strut 

(1.7) $\ \ k_{n}\left( x_{i_{1}}^{e_{1}},...,x_{i_{n}}^{e_{n}}\right)
=\oplus _{j=1}^{N}\left(
k_{n}^{(j)}(a_{i_{1}:j}^{e_{i_{1}}},...,a_{i_{n}:j}^{e_{i_{n}}})\right) ,$

\strut 

where $x_{i_{k}}^{e_{i_{k}}}=\oplus _{j=1}^{N}\left(
a_{i_{k}:j}^{e_{i_{k}}}\right) $ in $(A,$ $E),$ for all $k$ $=$ $1,$ ..., $n,
$ and where $e_{i_{1}},$ ..., $e_{i_{n}}$ $\in $ $\{1,$ $*\}$ and $(i_{1},$
..., $i_{n})$ $\in $ $\{1,$ $2\}^{n},$ for all $n$ $\in $ $\Bbb{N}.$
Therefore, for the mixed $n$-tuple $(x_{i_{1}}^{e_{1}},$ ..., $%
x_{i_{n}}^{e_{i_{n}}})$ of $\{x_{1},$ $x_{1}^{*},$ $x_{2},$ $x_{2}^{*}\},$
the $n$-tuples $(a_{i_{1}:k}^{e_{i_{1}}},$ ..., $a_{i_{n}:k}^{e_{i_{n}}})$
are also mixed $n$-tuple of $\{a_{k},$ $a_{k}^{*},$ $b_{k},$ $b_{k}^{*}\},$
for all $k$ $=$ $1,$ ..., $N.$ Therefore, for such mixed $n$-tuple,

\strut 

(1.8) \ \ \ \ \ \ \ \ \ \ $k_{n}^{(k)}\left(
a_{i_{1}:k}^{e_{i_{1}}},...,a_{i_{n}:k}^{e_{i_{n}}}\right) =0,$ for all $k$ $%
=$ $1,...,N.$

\strut 

By (1.8), the $n$-th mixed trivial $D_{N}$-valued cumulants of $x_{1},$ $%
x_{1}^{*},$ $x_{2}$ and $x_{2}^{*}$ in (1.7) vanish, and hence the $D_{N}$%
-valued random variables $x_{1}$ and $x_{2}$ are free over $D_{N}$ in the
direct producted $*$-probability space $(A,$ $E).$
\end{remark}

\strut \strut

The above theorem and remark shows that the $D_{N}$-freeness of $\oplus
_{j=1}^{N}$ $a_{j}$ and $\oplus _{j=1}^{N}$ $b_{j}$ in the direct producted
noncommutative (or $*$- \ or $\ C^{*}$- \ or $\ W^{*}$-) probability space $%
(A,$ $E)$ is characterized by the (scalar-valued) freeness of $a_{j}$ and $%
b_{j}$ in $(A_{j},$ $\varphi _{j}),$ for all $j$ $=$ $1,$ ..., $N.$

\strut

\begin{corollary}
Let $e_{i}=0\oplus ...\oplus 0\oplus a_{i}\oplus 0\oplus ...\oplus 0$ and $%
e_{j}=0\oplus $ $...$ $\oplus $ $0$ $\oplus $ $a_{j}$ $\oplus $ $0$ $\oplus $
$...$ $\oplus $ $0$ in $(A,$ $E),$ where $a_{k}$ $\in $ $(A_{k},$ $\varphi
_{k}),$ for $k$ $=$ $i,$ $j.$ If $i$ $\neq $ $j,$ then $e_{i}$ and $e_{j}$
are free over $D_{N}$ in $(A,$ $E).$\ $\square $
\end{corollary}

\strut

\strut Define subalgebras $A_{1}^{\prime },$ $...,$ $A_{N}^{\prime }$ of the
direct product $A=\oplus _{j=1}^{N}A_{j}$ by

\strut

\begin{center}
$A_{j}^{\prime }=0\oplus ...\oplus 0\oplus A_{j}\oplus 0\oplus ...\oplus 0,$
\end{center}

\strut

for all $j$ $=1,...,N.$ Then $A_{j}^{\prime }$ is the embedding of $A_{j}$
in $A.$ By the previous corollary, we can easily get the following
proposition;

\strut

\begin{corollary}
The unital\ algebras $A_{1},$ ..., $A_{N}$ are free over $D_{N}$ in the
direct producted noncommutative probability space $(A,$ $E).$ $\square $
\end{corollary}

\strut

In the above corollary, we can replace the condition [algebras $A_{1},$ ..., 
$A_{N}$] to [$*$-algebras $A_{1},$ ..., $A_{N}$]. By definition, the unital
algebra $A_{j}$ is always free from $\Bbb{C}$ (See [1], [4], [10] and [11]).

\strut \strut \strut

\strut

\strut

\section{Graph Free Probability Spaces}

\strut

\strut

In this chapter, we will consider our main objects of this paper. Throughout
this chapter, let $G$ be a finite simplicial graph with its finite vertex
set $V(G)$ and the edge set $E(G).$ Let $e$ $\in $ $E(G)$ be an edge
connecting the vertices $v_{1}$ and $v_{2}.$ Then denote $e$ by $[v_{1},$ $%
v_{2}].$ Assume that $[v_{1},$ $v_{2}]$ $=$ $[v_{2},$ $v_{1}],$ if there
exists an edge $e$ connecting $v_{1}$ and $v_{2}.$ Since the graph $G$ is
simplicial, if $[v_{1},$ $v_{2}]$ is in $E(G),$ then this is the unique edge
connecting the vertices $v_{1}$ and $v_{2}.$ By $FP(G),$ we will denote the
set of all admissible finite paths. Then this set $FP(G)$ is partitioned by

\strut

\begin{center}
$FP(G)=\cup _{n=1}^{\infty }FP_{n}(G),$
\end{center}

with

\begin{center}
$FP_{n}(G)=\{w\in FP(G):\left| w\right| =n\},$
\end{center}

\strut

where $\left| w\right| $ is the length of the finite path $w$ in $FP(G).$
i.e., if $w$ $=e_{1}$ $e_{2}$ ... $e_{k},$ where $e_{1}$ $=$ $[v_{1},$ $%
v_{2}],$ $e_{2}$ $=$ $[v_{2},$ $v_{3}],$ ..., $e_{k}$ $=$ $[v_{k},$ $%
v_{k+1}] $ are admissible edges making the finite path $w,$ then the path $w$
is denoted by $[v_{1},$ $v_{2},$ ..., $v_{k+1}],$ and $\left| w\right| $ is
defined to be $k.$ It means that the finite path $w$ connects vertices $%
v_{1} $ and $v_{k+1}$ via $v_{2},$ ..., $v_{k-1},$ and the length $\left|
w\right| $ of $w$ is $k$ $=$ $\left| \{v_{1},\text{ ..., }v_{k+1}\}\right| $ 
$-$ $1.$ And then $w$ $\in $ $FP_{k}(G)$ in $FP(G).$ Clearly, the edge set $%
E(G)$ is $FP_{1}(G).$ The following is automatically assumed;

\strut

\begin{center}
$[v_{1},v_{2},...,v_{n-1},v_{n}]=[v_{n},v_{n-1},...,v_{2},v_{1}],$
\end{center}

\strut

for all $n\in \Bbb{N}$ $\setminus $ $\{1\},$ where $v_{1},$ ..., $v_{n}$ $%
\in $ $V(G).$ Define the free semigroupoid $\Bbb{F}^{+}(G)$ by

\strut

\begin{center}
$\Bbb{F}^{+}(G)=V(G)\cup FP(G).$
\end{center}

\strut

\strut

\strut

\subsection{Graph Free Product of Probability-Space-Vertices}

\strut

\strut

Throughout this section, let $G$ be a finite simplicial graph, having its
vertex-set $V(G)$ $=$ $\{1,$ ..., $N\},$ with the probability-space vertices 
$(A_{1},$ $\varphi _{1}),$ ..., $(A_{N},$ $\varphi _{N}).$ i.e., there
exists a set of noncommutative probability spaces, indexed by the given
vertex-set $V(G),$

\strut

\begin{center}
$\{(A_{1},\varphi _{1}),$ ..., $(A_{N},$ $\varphi _{N})\}.$
\end{center}

\strut

Such graphs with probability-space vertices are said to be graphs of
probability-space vertices.

\strut

\begin{definition}
Let $G$ be a graph with probability-space vertices $\{(A_{v},$ $\varphi
_{v}):$ $v$ $\in $ $V(G)\}.$ Define the graph free product $A^{G}$ of $%
\{A_{v}\}_{v\in V(G)},$ by the algebra

\strut 

$\ \ \ \ \ \ \ \ \ \ \ \ \ \ \ \ \ A^{G}=\Bbb{C}\oplus \left( \underset{w\in %
\Bbb{F}^{+}(G)}{\oplus }A_{w}\right) ,$

\strut 

where $A_{w}=A_{w},$ for all $w\in V(G),$ and

\strut 

$\ \ \ \ \ \ \ \ \ \ \ \ \ \ \ \ A_{w}=A_{v_{1}}*A_{v_{2}}*...*A_{v_{k}},$

\strut 

for all $w=[v_{1},v_{2},...,v_{k}]\in FP(G),$ $k\in \Bbb{N}.$ If $A_{w}$ is
a summand of $A^{G},$ we write $A_{w}$ $<_{G}$ $A^{G}.$ Sometimes, we will
denote $A^{G}$ by $\underset{v\in V(G)}{*^{G}}$ $A_{v}.$ The symbol ``$*^{G}$%
'' is called the graph free product. Since we have the free product $\varphi
_{w}$ $=$ $*_{n=1}^{k}$ $\varphi _{v_{n}}$ of linear functionals $\varphi
_{v_{1}},$ ..., $\varphi _{v_{n}},$ for each finite path $w$ $=$ $[v_{1},$
..., $v_{k}],$ we have the corresponding noncommutative probability space $%
(A_{w},$ $\varphi _{w}).$ (Here, the symbol ``$*$'' means the usual free
product.) Now, define the $G$-diagonal subalgebra $D^{G}$ of $A^{G}$ by

\strut 

$\ \ \ \ \ \ \ \ \ \ \ \ \ \ \ \ \ D^{G}=\Bbb{C}\oplus \left( \underset{w\in %
\Bbb{F}^{+}(G)}{\oplus }\Bbb{C}_{w}\right) ,$

\strut 

where $\Bbb{C}_{v}=\Bbb{C}=\Bbb{C}_{w},$ for all $v\in V(G)$ and $w\in FP(G).
$ Then, like in Chapter 1, we can define the conditional expectation $E^{G}$
from $A^{G}$ into $D^{G}$ by

\strut 

$\ \ \ \ \ \ \ \ \ \ \ \ \ \ \ \ \ \ \ \ \ E^{G}=\left( \underset{w\in %
\Bbb{F}^{+}(G)}{\oplus }\varphi _{w}\right) .$

\strut 

Then the algebraic pair $\left( A^{G},E^{G}\right) $ is a noncommutative
probability space with amalgamation over $D^{G},$ and it is called the 
\textbf{graph free probability space} of $\{(A_{v},$ $\varphi _{v})$ $:$ $v$ 
$\in $ $V(G)\}$. All elements in the graph free probability space $(A^{G},$ $%
E^{G})$ are called graph random variables or $G$-random variables.
\end{definition}

\strut

By definition, we can easily check that the subalgebra $D^{G}$ commutes with 
$A^{G}.$ i.e.,

\strut

(2.1) \ \ \ \ \ \ \ $\ \ da=ad,$ for all $d\in D^{G}$ and $a\in A^{G}.$

\strut

By definition and (2.1), we have that;

\strut

\begin{proposition}
Let $(A^{G},E^{G})$ be a graph free probability space over its subalgebra $%
D^{G}.$ Then it is a direct producted noncommutative probability space of $%
\{(A_{w},$ $\varphi _{w})$ $:$ $w$ $\in $ $\Bbb{F}^{+}(G)\}.$ $\square $
\end{proposition}

\strut

\begin{definition}
(\textbf{Graph-Freeness}) Let $G$ be a graph with probability-space vertices
and let $(A^{G},$ $E^{G})$ be the corresponding graph free probability space
over its $G$-diagonal subalgebra $D^{G}.$ The subalgebras $A_{1}$ and $A_{2}$
of $A^{G}$ are said to be graph-free or $G$-free if all mixed $D^{G}$-valued
cumulants of $A_{1}$ and $A_{2}$ vanishes. The subsets $X_{1}$ and $X_{2}$
of $A^{G}$ are $G$-free if the subalgebras Alg $(X_{1},$ $D^{G})$ and Alg $%
(X_{2},$ $D^{G})$ are $G$-free.
\end{definition}

\strut

By the previous definition, we have that;

\strut

\begin{theorem}
Let $x$ and $y$ be $G$-random variables in $(A^{G},$ $E^{G}).$ The $G$%
-random variables $x$ and $y$ are $G$-free if and only if either (i) or (ii)
holds;

\strut \strut 

(i) \ there exists $(A_{w},$ $\varphi _{w})<_{G}(A^{G},E^{G})$ such that $x$
and $y$ are free in $(A_{w},$ $\varphi _{w}).$

(ii) they are free over $D^{G}$ in $(A^{G},$ $E^{G}).$
\end{theorem}

\strut

\begin{proof}
Let $x$ and $y$ be $G$-random variables in the graph free probability space $%
(A^{G},$ $E^{G}).$ Then, by the definition of $A_{G},$

\strut

\ \ \ \ \ \ \ $\ \ \ \ \ \ x=\underset{w\in \Bbb{F}^{+}(G)}{\oplus }x_{w}$ \
\ \ and \ \ \ $y=\underset{w^{\prime }\in \Bbb{F}^{+}(G)}{\oplus }%
y_{w^{\prime }}.$

\strut

($\Rightarrow $) Assume that the $G$-random variables $x$ and $y$ are $G$%
-free. i.e., they have the vanishing mixed $D^{G}$-valued cumulants. By
Chapter 1, all summands $x_{w}$ and $y_{w^{\prime }}$ of $x$ and $y$ are
free in $(A_{w},$ $\varphi _{w}).$ So, if there exists $w_{0}$ $\in $ $\Bbb{F%
}^{+}(G)$ such that $x$ $=$ $x_{w_{0}}$ and $y$ $=$ $y_{w_{0}},$ then, as
scalar-valued random variables, $x$ and $y$ are free in $(A_{w_{0}},$ $%
\varphi _{w_{0}})$ $<_{G}$ $(A^{G},$ $E^{G}).$ Otherwise, again by Chapter
1, they are free over $D^{G}$ in $(A^{G},$ $E^{G}),$ as $D^{G}$-valued
random variables in the direct producted noncommutative probability space $%
A^{G}$ of $A_{w}$'s, for $w$ $\in $ $\Bbb{F}^{+}(G).$

($\Leftarrow $) Suppose the $G$-random variables $x$ and $y$ are free in $%
(A_{w},$ $\varphi _{w})$, where $A_{w}$ $<_{G}$ $A^{G}.$ First, this means
that the operators $x$ and $y$ are contained in the summand $A_{w}$ of $%
A^{G}.$ Also, since they are free in $(A_{w},$ $\varphi _{w}),$ they are
free over $D^{G}$ in $(A^{G},$ $E^{G}),$ as $D^{G}$-valued random variables $%
x$ $=$ $x_{w}$ and $y$ $=$ $y_{w}.$ Thus they are $G$-free. Otherwise, if $x$
and $y$ are free over $D^{G}$ in $(A^{G},$ $E^{G}),$ then, by the very
definition of $G$-freeness, they are $G$-free in $(A^{G},$ $D^{G}).$\strut
\end{proof}

\strut \strut \strut \strut 

\strut \strut \strut \strut

\strut

\subsection{Graph-Freeness}

\strut

\strut

In this section, we will consider the graph-freeness more in detail.
Throughout this section, let $G$ be a graph with probability-space vertices

\strut

\begin{center}
$\{(A_{v},\varphi _{v}):v\in V(G)\}.$
\end{center}

\strut

Also, let $(A^{G},$ $E^{G})$ be the corresponding graph free probability
space over its $G$-diagonal subalgebra $D^{G}.$

\strut

\begin{lemma}
Let $A_{v_{1}}$ and $A_{v_{2}}$ be the direct summands of $A^{G},$ where $%
v_{1}$ $\neq $ $v_{2}$ in $V(G).$ Then they are $G$-free in $(A^{G},$ $%
E^{G}).$
\end{lemma}

\strut

\begin{proof}
In Chapter 1, we showed that the direct summands of a direct producted
noncommutative probability space are free from each other over the diagonal
subalgebra. So, if $v_{1}$ $\neq $ $v_{2}$ in $V(G),$ then $A_{v_{1}}$ and $%
A_{v_{2}}$ are free over $D^{G}$ in $(A^{G},$ $E^{G}),$ as direct summands
of $A^{G}.$ Equivalently, they are $G$-free in $(A^{G},$ $E^{G}).$
\end{proof}

\strut

We will consider more general case. To do that, we need the following
concept;

\strut

\begin{definition}
Let $\Bbb{F}^{+}(G)$ be the free semigroupoid of the given graph $G$ and let 
$w_{1}$ and $w_{2}$ be elements in $\Bbb{F}^{+}(G).$ We say that $w_{1}$ $=$ 
$[v_{1},$ ..., $v_{k}]$ and $w_{2}$ $=$ $[v_{1}^{\prime },$ ..., $%
v_{l}^{\prime }]$ are disjoint if $\{v_{1},$ ..., $v_{k}\}\cap
\{v_{1}^{\prime },...,v_{l}^{\prime }\}=\emptyset ,$ for $k,$ $l\in \Bbb{N}.$
(Notice that if $k$ $=$ $1$ and $l$ $=$ $1,$ then $w_{1}$ $=$ $[v_{1}]$ and $%
w_{2}$ $=$ $[v_{1}^{\prime }]$ are vertices in $V(G).$)
\end{definition}

\strut

\begin{theorem}
Let $B_{w_{1}}$ $\simeq $ $A_{w_{1}}$ and $B_{w_{2}}$ $\simeq $ $A_{w_{2}}$
be subalgebras of $A^{G}$ generated by $w_{1}$ and $w_{2},$ respectively.
(Notice that $A_{w_{1}}$ and $A_{w_{2}}$ are direct sums of $A^{G},$ but $%
B_{w_{1}}$ and $B_{w_{2}}$ are not necessarily direct sums. They are just
subalgebras of $A^{G}.$) The subalgebras $B_{w_{1}}$ and $B_{w_{2}}$ are $G$%
-free if and only if either (i) $B_{w_{1}}$ $=$ $A_{w_{1}}$ and $B_{w_{2}}$ $%
=$ $A_{w_{2}}$ or (ii) $w_{1}$ and $w_{2}$ are disjoint.
\end{theorem}

\strut

\begin{proof}
($\Leftarrow $) If $B_{w_{1}}$ $=$ $A_{w_{1}}$ and $B_{w_{2}}$ $=$ $%
A_{w_{2}},$ then they are $G$-free, by definition. If $%
w_{1}=[v_{1},...,v_{k}]$ and $w_{2}=[v_{1}^{\prime },...,v_{l}^{\prime }]$
are disjoint ($k,$ $l$ $\in $ $\Bbb{N}$), then the direct summands $A_{w_{1}}
$ and $A_{w_{2}}$ of the graph free product $A^{G}$ have no common direct
summands. This means that, for any pair $(v_{i},$ $v_{j}^{\prime })$ $\in $ $%
\{v_{1},$ ..., $v_{k}\}$ $\times $ $\{v_{1}^{\prime },$ ..., $v_{l}^{\prime
}\},$ the direct summands $A_{v_{i}}$ and $A_{v_{j}^{\prime }}$ are $G$%
-free, by the previous lemma. Therefore,

\strut

$\ \ \ \ \ \ \ \ \ \ B_{w_{1}}=A_{v_{1}}*...*A_{v_{k}}$ \ and \ $%
B_{w_{2}}=A_{v_{1}^{\prime }}*...*A_{v_{l}^{\prime }}$

\strut

are free over $D^{G}$ in $(A^{G},$ $E^{G}),$ and hence they are $G$-free.

\strut

($\Rightarrow $) Suppose $B_{w_{1}}$ and $B_{w_{2}}$ are $G$-free in $(A^{G},
$ $E^{G})$ and assume that $w_{1}$ $=$ $[v_{1},$ ..., $v_{k}]$ and $w_{2}$ $=
$ $[v_{1}^{\prime },$ ..., $v_{l}^{\prime }]$ are not disjoint in $\Bbb{F}%
^{+}(G),$ where $k,$ $l$ $\in $ $\Bbb{N}.$ Then, by definition, there exists
the nonempty intersection of $\{v_{1},$ ..., $v_{k}\}$ and $\{v_{1}^{\prime
},$ ..., $v_{l}^{\prime }\}.$ Take a vertex $v_{0}$ in the intersection.
Then the algebra $A_{v_{0}}$ of $A^{G}$ is contained in both $A_{w_{1}}$ and 
$A_{w_{2}}.$ i.e.,

\strut

\ \ \  $\ \ \ \ \ \ \ \ \ \ \ \ \ \ \ \ \ B_{w_{1}}=A_{v_{0}}*\left( 
\underset{v\in \{v_{1},...,v_{k}\}\,\setminus \,\{v_{0}\}}{*}A_{v}\right) $

and

\ \ \  $\ \ \ \ \ \ \ \ \ \ \ \ \ \ \ \ B_{w_{2}}=A_{v_{0}}*\left( \underset{%
v^{\prime }\in \{v_{1}^{\prime },...,v_{l}^{\prime }\}\,\setminus \,\{v_{0}\}%
}{*}A_{v^{\prime }}\right) .$

\strut

(Recall that $A*B=B*A,$ for algebras $A$ and $B.$) We can take $x_{0}$ $\in $
$A_{v_{0}}$ $\setminus $ $\Bbb{C}$ such that $\varphi _{v_{0}}(x_{0}^{k})$ $%
\neq $ $0,$ for some $k$ $\in $ $\Bbb{N}.$ Define $x$ $=$ $x_{0}$ $\in $ $%
B_{w_{1}}$ and $y$ $=$ $x_{0}$ $\in $ $B_{w_{2}}.$ Then the random variables 
$x$ and $y$ have the nonvanising mixed cumulants, with respect to the linear
functional $\varphi _{v_{0}}.$ Therefore, by regarding them as $G$-random
variables, they have the nonvanishing mixed $D^{G}$-valued cumulants, with
respect to the conditional expectation $E^{G}$. This shows that $B_{w_{1}}$
and $B_{w_{2}}$ are not free over $D^{G}$ in the direct producted
noncommutative probability space $(A^{G},$ $E^{G}).$ Equivalently, they are
not $G$-free. This contradict our assumption that $B_{w_{1}}$ and $B_{w_{2}}$
are $G$-free.
\end{proof}

\strut \strut

By the previous theorem, we can get the following corollaries;

\strut

\begin{corollary}
Let $W_{1}$ and $W_{2}$ are subsets in $\Bbb{F}^{+}(G)$ and assume that $%
W_{1}$ and $W_{2}$ are disjoint, in the sense that all pairs $(w_{1},$ $%
w_{2})$ $\in $ $W_{1}$ $\times $ $W_{2}$ are disjoint. Then the subalgebras $%
\underset{w_{1}\in W_{1}}{\oplus }$ $A_{w_{1}}$ and $\underset{w_{2}\in W_{2}%
}{\oplus }$ $A_{w_{2}}$ are $G$-free in $(A^{G},$ $E^{G}).$ And the converse
also holds true. $\square $
\end{corollary}

\strut

\begin{corollary}
Let $w=[v_{1},...,v_{j},...,v_{k}]$ be a finite path, for $1\leq j<k$ in $%
\Bbb{N}.$ Then $A_{[v_{1},...,v_{j}]}$ and $A_{[v_{j+1},...,v_{k}]}$ are $G$%
-free in $(A^{G},$ $E^{G}).$
\end{corollary}

\strut

\begin{proof}
Notice that since $[v_{1},...,v_{k}]$ is a finite path, $[v_{1},...,v_{j}]$
and $[v_{j+1},$ $...,$ $v_{k}]$ are also in $\Bbb{F}^{+}(G).$ Since $[v_{1},$
..., $v_{j}]$ and $[v_{j+1},$ ..., $v_{k}]$ are disjoint, as direct summands
of $A^{G},$ the algebras $A_{[v_{1},...,v_{j}]}$ and $%
A_{[v_{j+1},...,v_{k}]} $ are $G$-free.
\end{proof}

\strut \strut \strut 

\strut \strut \strut \strut

\strut \strut \strut

\subsection{Graph R-transform Calculus}

\strut

\strut

In this section, we will define graph moment series and graph R-transforms
of graph random variables. Like before, let $G$ be a graph with
probability-space vertices and let $(A^{G},$ $E^{G})$ be the corresponding
graph free probability space over the $G$-diagonal subalgebra $D^{G},$ where 
$A^{G}$ $=$ $\underset{w\in \Bbb{F}^{+}(G)}{*^{G}}$ $A_{w}$ and $E^{G}$ $=$ $%
\underset{w\in \Bbb{F}^{+}(G)}{\oplus }$ $\varphi _{w}.$ 

\strut \strut 

\begin{definition}
Let $(A^{G},$ $E^{G})$ be a graph free probability space over the $G$%
-diagonal sublagebra $D^{G}$ and let $x$ $\in $ $(A^{G},$ $E^{G})$ be a $G$%
-random variable. The graph moments (or $G$-moments) of $x$ and the graph
cumulants (or $G$-cumulants) of $x$ are defined by

\strut 

$\ \ \ \ \ \ \ \ \ \ \ \ \ E^{G}(x^{n})$ \ \ and \ \ $k_{n}^{(E^{G})}\left( 
\underset{n\text{-times}}{\underbrace{x,.......,x}}\right) ,$

\strut 

for all $n\in \Bbb{N}.$\strut 
\end{definition}

\strut

Let $x$ be a $G$-random variable in $(A^{G},$ $E^{G}).$ Then $x$ $=$ $%
\underset{w\in \Bbb{F}^{+}(G)}{\oplus }$ $x_{w}.$ So, by Chapter 1, we have
that

\strut

(2.2) \ $\ \ \ \ \ \ \ \ \ \ \ \ \ \ \ \ \ \ E^{G}\left( x^{n}\right) =%
\underset{w\in \Bbb{F}^{+}(G)}{\oplus }\left( \varphi _{w}(x_{w}^{n})\right)
,$

\strut

and

\strut

(2.3) \ $\ \ \ \ \ \ \ k_{n}^{(E^{G})}\left( x,...,x\right) =\underset{w\in 
\Bbb{F}^{+}(G)}{\oplus }\left( k_{n}^{(\varphi _{w})}\left(
x_{w},...,x_{w}\right) \right) ,$

\strut

for all $x=\underset{w\in \Bbb{F}^{+}(G)}{\oplus }x_{w}\in A^{G}$ and for
all $n\in \Bbb{N}.$ Notice that, by (2.1), we have that

\strut

\begin{center}
$k_{n}^{(E^{G})}(x,...,x)=\underset{\pi \in NC(n)}{\sum }\left( \underset{%
V\in \pi }{\Pi }E^{G}(x^{\left| V\right| })\right) \mu (\pi ,1_{n}),$
\end{center}

\strut

like the scalar-valued cumulants. Recall that, when we consider the
operator-valued M\"{o}bius inversion, we have to think about the insertion
property. However, in our case, by (2.1), we need not think about the
insertion property. Also, the trivial $D^{G}$-valued cumulants (or moments)
of $x$ contains the full free distributional data of $x.$

\strut

\begin{definition}
Let $x\in (A^{G},$ $E^{G})$ be a $G$-random variable. Define the graph
moment series (or $G$-moment series) of $x$ by

\strut 

$\ \ \ \ \ \ \ \ \ \ \ \ \ \ \ \ \ \ \ M_{x}^{G}(z)=\sum_{n=1}^{\infty
}\left( E^{G}(x^{n})\right) z^{n}$

\strut 

in the ring $D^{G}[[z]]$ of formal series, where $z$ is an indeterminent..
Also, define the graph R-transform (or $G$-R-transform) of $x$ by

\strut 

$\ \ \ \ \ \ \ \ \ \ \ \ \ \ \ R_{x}^{G}(z)=\sum_{n=1}^{\infty }\left(
k_{n}^{(E^{G})}(x,...,x)\right) z^{n}$

\strut 

in $D^{G}[[z]].$ More generally, if $x_{1},$ ..., $x_{s}$ are $G$-random
variables, then the $G$-moment series of $x_{1},$ ..., $x_{s}$ and $G$%
-R-transform of them are defined by

\strut 

$\ \ \ M_{x_{1},...,x_{s}}^{G}(z_{1},...,z_{s})=\sum_{n=1}^{\infty }%
\underset{(i_{1},...,i_{n})\in \{1,...,s\}^{n}}{\sum }%
E^{G}(x_{i_{1}}...x_{i_{n}})\,z_{i_{1}}...z_{i_{n}}$

and

$\ \ \ R_{x_{1},...,x_{s}}^{G}(z_{1},...,z_{s})=\sum_{n=1}^{\infty }%
\underset{(i_{1},...,i_{n})\in \{1,...,s\}^{n}}{\sum }k_{n}^{(E_{G})}\left(
x_{i_{1}},...,x_{i_{n}}\right) z_{i_{1}}...z_{i_{n}},$

\strut 

in $D^{G}[[z_{1},...,z_{s}]],$ where $z_{1},$ ..., $z_{s}$ are
noncommutative indeterminants.
\end{definition}

\strut

By (2.1), the graph moment series and graph R-transforms of graph random
variables are well-defined. 

\strut

\begin{proposition}
Let $x$ and $y$ be $G$-random variables and assume that they are $G$-free.
Then

\strut 

(1) $\ \ R_{x+y}^{G}(z)=\left( R_{x}^{G}+R_{y}^{G}\right) (z)$

(2) \ \ $R_{x,y}^{G}(z_{1},z_{2})=R_{x}^{G}(z_{1})+R_{y}^{G}(z_{2}).$ \ $%
\square $
\end{proposition}

\strut

More generally, we have that;

\strut

\begin{theorem}
Let $\{x_{1},...,x_{s}\}$ and $\{y_{1},...,y_{s}\}$ be sets of $G$-random
variables ($s$ $\in $ $\Bbb{N}$) in the graph free probability space $(A^{G},
$ $E^{G}).$ If these two sets are $G$-free, then

\strut 

(1) $R_{x_{1}+y_{1},...,x_{s}+y_{s}}^{G}(z_{1},...,z_{s})=$ $\left(
R_{x_{1},...,x_{s}}^{G}+R_{y_{1},...,y_{s}}^{G}\right) (z_{1},...,z_{s})$

(2) $%
R_{x_{1},...,x_{s},y_{1},...,y_{s}}^{G}(z_{1},..,z_{2s})=R_{x_{1},...,x_{s}}^{G}(z_{1},...,z_{s})+R_{y_{1},...,y_{s}}^{G}(z_{s+1},...,z_{2s}).
$ $\square $
\end{theorem}

\strut

Let $x$ and $y$ be $G$-random variables in the graph free probability space $%
(A^{G},$ $E^{G})$ and assume that they are $G$-free. Then they are free over 
$D^{G}$ in $(A^{G},$ $E^{G}),$ by regarding $(A^{G},$ $E^{G})$ as a direct
producted noncommutative probability space. So, we have that

\strut

$\ \ \ \ \ \ \ k_{n}^{(E^{G})}\left( xy,...,xy\right) $

\strut

(2.4) $\ \ \ \ \ \ \ \ \ \ \ =\underset{\pi \in NC(n)}{\sum }\left( k_{\pi
}^{(E^{G})}(x,...,x)\right) \left( k_{Kr(\pi )}^{(E^{G})}(y,...,y)\right) ,$

\strut \strut \strut

for $n\in \Bbb{N},$ where $Kr:NC(n)\rightarrow NC(n)$ is the Kreweras
complementation map on $NC(n).$ In fact,

\strut 

\strut (2.5)

\begin{center}
$k_{n}^{(E^{G})}\left( xy,...,xy\right) =\underset{\pi \in NC(n)}{\sum }%
k_{\pi \cup _{alt}Kr(\pi )}^{(E_{G})}\left( x,y,x,y,...,x,y\right) ,$
\end{center}

\strut

by the $D^{G}$-freeness of $x$ and $y,$ where $\pi \cup _{alt}\theta $ means
the alternating sum of partitions $\pi $ and $\theta $ in $NC(n),$ $n$ $\in $
$\Bbb{N}.$ i.e., if

\strut

\begin{center}
$\pi =\{(1,4,5),(2,3),(6,8),(7)\}$ $\in NC(8),$\strut
\end{center}

then\strut

\begin{center}
$Kr(\pi )=\{(1,3),(2),(4),(5,8),(6,7)\}\in NC(8)$
\end{center}

and

\begin{center}
$
\begin{array}{ll}
\pi \cup _{alt}Kr(\pi )= & \{(1,7,9),(2,6),(3,5),(4), \\ 
& \text{ \ }(10,16),(11,15),(12,14),(13)\}\in NC(16).
\end{array}
$
\end{center}

\strut

However, by (2.1), the formula (2.5) is exactly same as (2.4). \strut So,
the $G$-R-transform $R_{xy}^{G}(z)$ of $xy$ satisfies that

\strut

(2.6)

\begin{center}
$R_{xy}^{G}(z)=\sum_{n=1}^{\infty }\left( \underset{\pi \in NC(n)}{\sum }%
\left( k_{\pi }^{(E^{G})}(x,...,x)\right) \left( k_{Kr(\pi
)}^{(E^{G})}(y,...,y)\right) \right) z^{n}.$
\end{center}

\strut

More generally, if two sets $\{x_{1},$ ..., $x_{s}\}$ and $\{y_{1},$ ..., $%
y_{s}\}$ are $G$-free in $(A^{G},$ $E^{G}),$ then

\strut

$k_{n}^{(E^{G})}\left(
x_{i_{1}}y_{i_{1}},x_{i_{2}}y_{i_{2}},...,x_{i_{n}}y_{i_{n}}\right) $

\strut

(2.7) \ \ \ \ \ \ \ $=\underset{\pi \in NC(n)}{\sum }\left( k_{\pi
}^{(E^{G})}(x_{i_{1}},...,x_{i_{n}})\right) \left( k_{Kr(\pi
)}^{(E^{G})}(y_{i_{1}},...,y_{i_{n}})\right) $,

\strut

where $(i_{1},...,i_{n})\in \{1,...,s\}^{n},$ for $n\in \Bbb{N}.$ So,

\strut

(2.8) $\ \ \ R_{x_{1}y_{1},...,x_{s}y_{s}}^{G}(z_{1},...,z_{s})=$

$\strut $

\begin{center}
$\sum_{n=1}^{\infty }\underset{(i_{1},...,i_{n})\in \{1,...,s\}^{n}}{\sum }%
\left( \underset{\pi \in NC(n)}{\sum }\left( k_{\pi
}^{(E^{G})}(x_{i_{1}},...,x_{i_{n}})\right) \left( k_{Kr(\pi
)}^{(E^{G})}(y_{i_{1}},...,y_{i_{n}})\right) \right) z_{i_{1}}...z_{i_{n}}$
\end{center}

\strut \strut

\begin{quote}
\textbf{Notation} \ In the rest of this paper, the right-hand side of (2.6)
is denoted by

\strut
\end{quote}

\begin{center}
$\left( R_{x}^{G}\,\,\,\,\frame{*}_{G}\,\,\,R_{y}^{G}\right) (z).$
\end{center}

\begin{quote}
\strut

And similarly, the right-hand side of (2.8) is denoted by

\strut
\end{quote}

\begin{center}
$\left( R_{x_{1},...,x_{s}}^{G}\,\,\frame{*}_{G}\,\,R_{y_{1},...,y_{s}}^{G}%
\right) (z_{1},...,z_{s}).$
\end{center}

\begin{quote}
\strut

The symbol ``\frame{*}$_{G}$'' is motivated by the boxed convolution \frame{*%
}$_{s}$ of Nica. Recall that the Nica's boxed convolution \frame{*}$_{s}$
makes a certain subset $\Theta _{s}$ of $\Bbb{C}[[z_{1},$ ..., $z_{s}]]$ be
semigroup (See [1]). But our boxed-convolution-like notation \frame{*}$_{G}$
does not have such meaning. This is just a notation representing the
right-hand sides of (2.6) and (2.8). $\square $
\end{quote}

\strut

By the previous discussion, we have that;

\strut

\begin{proposition}
Let $\{x_{1},...,x_{s}\}$ and $\{y_{1},...,y_{s}\}$ be $G$-free in the $G$%
-free probability space $(A^{G},$ $E^{G}).$ Then

\strut \strut 

$\ \ \ \ \ \ \ R_{x_{1}y_{1},...,x_{s}y_{s}}^{G}(z_{1},...,z_{s})=\left(
R_{x_{1},...,x_{s}}^{G}\,\,\frame{*}_{G}\,\,R_{y_{1},...,y_{s}}^{G}\right)
(z_{1},...,z_{s}),$

\strut 

in $D^{G}[[z_{1},...,z_{s}]].$ In particular, for each $j$ $\in $ $\{1,$
..., $s\},$

\strut 

$\ \ \ \ \ \ \ \ \ \ \ \ \ \ \ \ \ \ \ R_{x_{j}y_{j}}^{G}(z)=\left(
R_{x}^{G}\,\,\,\frame{*}_{G}\,\,\,R_{y}^{G}\right) (z)$

\strut 

in $D^{G}[[z]].$ \ $\square $
\end{proposition}

\strut

By the previous propositions, we have the following graph R-transform
calculus; \ if $\{x_{1},$ ..., $x_{s}\}$ and $\{y_{1},$ ..., $y_{s}\}$ are $G
$-free, then

\strut

(1) $R_{x_{1}+y_{1},...,x_{s}+y_{s}}^{G}(z_{1},...,z_{s})=\left(
R_{x_{1},...,x_{s}}^{G}+R_{y_{1},...,y_{s}}^{G}\right) (z_{1},...,z_{s})$

(2) $%
R_{x_{1},...,x_{s},y_{1},...,y_{s}}^{G}(z_{1},...,z_{2s})=R_{x_{1},...,x_{s}}^{G}(z_{1},...,z_{s})+R_{y_{1},...,y_{s}}^{G}(z_{s+1},...,z_{2s}) 
$

(3) $R_{x_{1}y_{1},...,x_{s}y_{s}}^{G}(z_{1},...,z_{s})=\left(
R_{x_{1},...,x_{s}}^{G}\,\,\,\frame{*}_{G}\,\,\,R_{y_{1},...,y_{s}}^{G}%
\right) (z_{1},...,z_{s}),$

\strut

in $D^{G}[[z_{1},...,z_{s}]].$\ \ \ 

\strut \strut

\strut

\strut \strut

\section{Graph-Free Random Variables}

\strut

\strut

In this chapter, we will let each probability-space vertex be a $W^{*}$%
-probability space. Let $G$ be a graph with $W^{*}$-probability-space
vertices. i.e.,

\strut

\begin{center}
$\{(A_{v},\varphi _{v}):v\in V(G),$ $\varphi _{v}$ is a state$\}$
\end{center}

\strut

is a family of $W^{*}$-probability spaces. Define graph free product $A^{G}$
of $A_{v}$'s by

\strut

\begin{center}
$A^{G}=\Bbb{C}\overline{\oplus }\left( \underset{w\in \Bbb{F}^{+}(G)}{%
\overline{\oplus }}A_{w}\right) ,$
\end{center}

\strut

where $\overline{\oplus }$ is the corresponding topological direct sum.
Also, define the $G$-diagonal algebra $D^{G}$ by

\strut

\begin{center}
$D^{G}=\Bbb{C}\overline{\oplus }\left( \underset{w\in \Bbb{F}^{+}(G)}{%
\overline{\oplus }}\Bbb{C}_{w}\right) ,$ \ with \ $\Bbb{C}_{w}=\Bbb{C},$ $%
\forall w.$
\end{center}

\strut

Also, the conditional expectation is defined by $E^{G}=\underset{w\in \Bbb{F}%
^{+}(G)}{\overline{\oplus }}\varphi _{w}.$

\strut

\begin{quote}
\textbf{Notation} \ For convenience, if there is no confusion, we will keep
using the notation $\oplus $ instead of $\overline{\oplus }.$ $\square $
\end{quote}

\strut

Notice that the corresponding graph free probability space $(A^{G},$ $E^{G})$
is a $W^{*}$-probability space with amalgamation over the $G$-diagonal
algebra $D^{G}.$ In such amalgamated $W^{*}$-probability space $(A^{G},$ $%
E^{G}),$ we will consider certain $G$-random variables. Notice that if $x$
is a $G$-random variables in a graph free $W^{*}$-probability space $(A^{G},$
$E^{G}),$ then 

\strut

\begin{center}
$x=\underset{w\in \Bbb{F}^{+}(G)}{\oplus }x_{w}.$
\end{center}

\strut

i.e., we can understand $x$ is the (infinite) direct sum of $x_{w}$'s, where 
$w$ $\in $ $\Bbb{F}^{+}(G).$

\strut \strut 

\strut \strut

\subsection{Graph Semicircular Elements}

\strut

\strut

In this section, we will consider the graph-semicircularity of graph random
variables. Let $G$ be a finite simplicial graph with $W^{*}$%
-probability-space-vertices $\{(A_{v},$ $\varphi _{v})$ $:$ $v$ $\in $ $%
V(G)\}$ and let $(A^{G},$ $E^{G})$ be the corresponding graph free $W^{*}$%
-probability space over its $G$-diagonal subalgebra $D^{G},$ with the $G$%
-conditional expectation $E^{G}$ $=$ $\underset{w\in \Bbb{F}^{+}(G)}{\oplus }
$ $\varphi _{w}.$

\strut

\begin{definition}
Let $x\in (A^{G},$ $E^{G})$ be a self-adjoint $G$-random variable. We say
that this $G$-random variable $x$ is graph-semicircular (or $G$%
-semicircular) if the only second $G$-cumulant of $x$ is nonvanishing. i.e., 
$x$ has the following $G$-cumulant relation;

\strut 

$\ \ \ \ \ \ \ \ \ k_{n}^{(E^{G})}\left( \underset{n\text{-times}}{%
\underbrace{x,.......,x}}\right) =\left\{ 
\begin{array}{lll}
k_{2}^{(E^{G})}(x,x) &  & \text{if }n=2 \\ 
&  &  \\ 
0_{D^{G}} &  & \text{otherwise.}
\end{array}
\right. $
\end{definition}

\strut

Equivalently, the $G$-random variable $x$ is $G$-semicircular if $x$ is $%
D^{G}$-valued semicircular, as a $D^{G}$-valued random variable in the
direct producted $W^{*}$-probability space $(A^{G},$ $E^{G}).$ In the rest
of this section, we will observe the conditions when the $G$-random variable 
$x$ is $G$-semicircular.

\strut

\begin{lemma}
Let $(A_{v},$ $\varphi _{v})$ be a $W^{*}$-probability-space-vertex of the
graph $G.$ If $a$ $\in $ $(A_{v},$ $\varphi _{v})$ is (scalar-valued)
semicircular, then, as a $G$-random variable $a$ is $G$-semicircular.
\end{lemma}

\strut

\begin{proof}
In general, if $x=\underset{w\in \Bbb{F}^{+}(G)}{\oplus }x_{w}$ is a $G$%
-random variable in $(A^{G},$ $E^{G}),$ then

\strut

$\ \ \ \ \ \ \ \ \ k_{n}^{(E^{G})}\left( \underset{n\text{-times}}{%
\underbrace{x,......,x}}\right) =\underset{w\in \Bbb{F}^{+}(G)}{\oplus }%
\left( k_{n}^{(\varphi _{w})}(x_{w},...,x_{w})\right) ,$

\strut

by Chapter 1 and (2.3). Thus, if $a$ $\in (A_{v},\varphi _{v})<_{G}(A^{G},$ $%
E^{G}),$ then

\strut

(2.9)$\ \ \ \ \ \ \ \ \ \ \ \ \ \ k_{n}^{(E^{G})}\left( a,...,a\right)
=k_{n}^{(\varphi _{v})}(a,...,a).$

\strut

Since $a$ is semicircular in $(A_{v},$ $\varphi _{v}),$ as a $G$-random
variable, it is also $G$-semicircular, by (2.9).
\end{proof}

\strut

By the previous lemma, we can get that;

\strut

\begin{proposition}
Let $w=[v_{1},...,v_{k}]\in FP(G)$ and let nonzero $a_{v_{j}}\in
(A_{v_{j}},\varphi _{v_{j}})$ be semicircular, for $j$ $=$ $1,$ ..., $k.$
Then the $G$-random variable $\sum_{j=1}^{k}$ $a_{v_{j}}$ is $G$%
-semicircular in $(A^{G},$ $E^{G}).$
\end{proposition}

\strut

\begin{proof}
By assumption, if the random variables $a_{v_{j}}$'s are nonzero, then they
are semicircular in $(A_{v_{j}},$ $\varphi _{v_{j}}),$ for some $j$ $=$ $1,$
..., $k.$ Notice that $B_{v_{1}},$ ..., $B_{v_{k}}$ are $G$-free in $(A^{G},$
$E^{G}),$ where $B_{v_{j}}$'s are subalgebra isomorphic to the direct sums $%
A_{v_{j}},$ for all $j$ $=$ $1,$ ..., $k.$ Indeed, since $\{v_{1}\},$ ..., $%
\{v_{k}\}$ are mutually disjoint, the subalgebras $B_{v_{1}},$ ..., $%
B_{v_{k}}$ are $G$-free from each other. Therefore, as $G$-random variables $%
a_{v_{1}},$ ..., $a_{v_{k}}$ are $G$-free from each other. So,

\strut

$\ \ \ \ \ \ k_{n}^{(E^{G})}\left(
\sum_{j=1}^{N}a_{v_{j}},...,\sum_{j=1}^{N}a_{v_{j}}\right) $

\strut

$\ \ \ \ \ \ \ \ \ \ \ \ \ \ \ =\sum_{j=1}^{N}k_{n}^{(E^{G})}\left(
a_{v_{j}},...,a_{v_{j}}\right) $

\strut

$\ \ \ \ \ \ \ \ \ \ \ \ \ \ \ =\left\{ 
\begin{array}{lll}
\sum_{j=1}^{N}k_{2}^{\varphi _{v_{j}}}\left( a_{v_{j}},a_{v_{j}}\right)  & 
& \text{if }n=2 \\ 
&  &  \\ 
0_{D^{G}} &  & \text{otherwise}
\end{array}
\right. $

\strut

by the previous lemma

\strut

$\ \ \ \ \ \ \ \ \ \ \ \ \ \ \ =\left\{ 
\begin{array}{ll}
k_{2}^{(E^{G})}\left(
\sum_{j=1}^{N}a_{v_{j}},\,\,\sum_{j=1}^{N}a_{v_{j}}\right)  & \text{if }n=2
\\ 
&  \\ 
0_{D^{G}} & \text{otherwise,}
\end{array}
\right. $

\strut

for all $n\in \Bbb{N}.$
\end{proof}

\strut

\strut

\strut

\subsection{Graph Circular Elements}

\strut

\strut

In this section, we will consider the graph circularity on a graph free
probability space $(A^{G},$ $E^{G}).$

\strut

\begin{definition}
Let $a\in (A^{G},$ $E^{G})$ be a $G$-random variable and assume that there
exist self-adjoint $G$-random variables $x_{1}$ and $x_{2}$ in $(A^{G},$ $%
E^{G})$ such that $a$ $=$ $x_{1}$ $+$ $ix_{2},$ and the $G$-random variables 
$x_{1}$ and $x_{2}$ are $G$-semicircular elements in $(A^{G},$ $E^{G}).$ We
say that the $G$-random variable $a$ is $G$-circular if $x_{1}$ and $x_{2}$
are $G$-free in $(A^{G},$ $E^{G}).$ We also say that the pair $(x_{1},$ $%
x_{2})$ is the $G$-semicircular pair of $a.$
\end{definition}

\strut

Suppose that $G$-circular element $a=x_{1}+ix_{2}$ in $(A^{G},$ $E^{G})$ has 
$G$-semicircular elements $x_{1}$ and $x_{2}$ in $(A_{w},$ $\varphi _{w}),$
for $w$ $\in $ $\Bbb{F}^{+}(G).$ Then the $G$-random variable $a$ is also
contained in $(A_{w},$ $\varphi _{w})$ and it is circular, in the sense of
Voiculescu in this $W^{*}$-probability space $(A_{w},$ $\varphi _{w}).$ Such 
$G$-random variable $a$ is said to be a $w$\textbf{-circular element} in $%
(A^{G},$ $E^{G}).$ By the $G$-freeness characterization, if $w_{1}$ and $%
w_{2}$ are disjoint in $\Bbb{F}^{+}(G)$, then $B_{w_{1}}$ $\simeq $ $%
A_{w_{1}}$ and $B_{w_{2}}$ $\simeq $ $A_{w_{2}}$ are $G$-free in $(A^{G},$ $%
E^{G}).$ Let $x_{1}$ and $x_{2}$ be semicircular in $(B_{w_{1}},$ $\varphi
_{w_{1}})$ and $(B_{w_{2}},$ $\varphi _{w_{2}}),$ respectively. In this
case, the $G$-random variable $a$ $=$ $x_{1}$ $+$ $ix_{2}$ is again $G$%
-circular in $(A^{G},$ $E^{G}).$ We say that such $G$-circular element $a$
is $(w_{1},$\textbf{\ }$w_{2})$\textbf{-circular} in $(A^{G},$ $E^{G}).$
Notice that if $a$ is $(w_{1},$ $w_{2})$-circular, then $w_{1}$ and $w_{2}$
are not admissible.

\strut \strut 

\strut \strut

\strut

\subsection{Graph R-diagonal Elements}

\strut

\strut

In this section, we define the graph R-diagonality on $(A^{G},$ $E^{G}).$

\strut

\begin{definition}
Let $a\in (A^{G},E^{G})$ be a $G$-random variable. It is said to be $G$%
-R-diagonal if it has the following $G$-cumulant relation; the only
nonvanishing $G$-cumulants of $a$ are either

\strut 

$\ \ \ k_{2n}^{(E^{G})}\left( a,a^{*},a,a^{*},...,a,a^{*}\right) $ \ or \ $%
k_{2n}^{(E^{G})}(a^{*},a,a^{*},a...,a^{*},a),$

\strut 

for all $n\in \Bbb{N}.$ We will call the above $G$-cumulants the alternating 
$G$-cumulants of $a.$
\end{definition}

\strut

\strut \strut \strut \strut

\strut \strut \strut \strut \strut \strut

\strut \textbf{References}

\bigskip

\strut

\begin{quote}
{\small [1]\ \ \ A. Nica, R-transform in Free Probability, IHP course note.}

{\small [2]\ \ \ A. Nica, R-transforms of Free Joint Distributions and
Non-crossing Partitions, J. of Func. Anal, 135 (1996), 271-296.\strut }

{\small [3] \ \ D. Shlyakhtenko, Notes on Free Probability Theory, (2005),
Letcture Note, arXiv:math.OA/0504063v1.}

{\small [4]\ \ \ D.Voiculescu, K. Dykemma and A. Nica, Free Random
Variables, CRM Monograph Series Vol 1 (1992).\strut }

{\small [5] \ F. Radulescu, Singularity of the Radial Subalgebra of }$%
L(F_{N})${\small \ and the Puk\'{a}nszky Invariant, Pacific J. of Math, vol.
151, No 2 (1991)\strut , 297-306.\strut \strut }

{\small [6] \ \ I. Cho, Toeplitz Noncommutative Probability Spaces over
Toeplitz Matricial Algebras}$,${\small \ (2002), Preprint.\strut }

{\small [7] \ \ I. Cho, The Moment Series of the Generating Operator of }$%
L(F_{2})*_{L(F_{1})}L(F_{2})${\small , (2003), Preprint.}

{\small [8] \ \ I. Cho, Random Variables in a Graph }$W^{*}${\small %
-Probability Space, (2004), Ph. D thesis, Univ. of Iowa. }

{\small [9]\ \ \ I. Cho, Direct Producted Noncommutative Probability Spaces,
(2005), Preprint.}

{\small [10]\ K. J. Horadam, The Word Problem and Related Results for Graph
Product Groups, Proc. AMS, vol. 82, No 2, (1981) 157-164.}

{\small [11]\ R. Speicher, Combinatorics of Free Probability Theory IHP
course note }

{\small [12] R. Speicher, Combinatorial Theory of the Free Product with
Amalgamation and Operator-Valued Free Probability Theory, AMS Mem, Vol 132 ,
Num 627 , (1998).}
\end{quote}

\strut

\end{document}